\documentclass[reqno,11pt]{amsart}
\usepackage{amsmath,amssymb,mathrsfs,amsthm,amsfonts}
\usepackage[inline]{enumitem} 
\usepackage[usenames,dvipsnames]{xcolor}
\usepackage{hyperref}
\usepackage{comment}
\hypersetup{%
	colorlinks=true, linkcolor=blue,
	citecolor=ForestGreen
}
\usepackage[paper=letterpaper,margin=1in]{geometry}
\usepackage{acronym}
\usepackage{dsfont}

\newcommand{\Dim}{\ensuremath{\mathrm{Dim}_{\mathds{H}}}}

\renewcommand{\P}{\mathds{P}}

\newcommand{\E}{\mathds{E}}
\newcommand{\R}{\mathds{R}}

\newcommand{\Z}{\mathds{Z}}
\newcommand{\N}{\mathds{N}}

\newcommand{\G}{\mathcal{G}}
\renewcommand{\L}{\mathcal{L}}

\newcommand{\I}{\mathcal{I}}
\renewcommand{\S}{\mathcal{S}}
\newcommand{\V}{\mathcal{V}}
\newcommand{\Q}{\mathds{Q}}

\newcommand{\e}{\mathrm{e}}

\newcommand{\deq}{\overset{\underset{\mathrm{d}}{}}{=}}

\newtheorem{stat}{Statement}[section]
\newtheorem{proposition}[stat]{Proposition}

\newtheorem{conjecture}[stat]{Conjecture}
\newtheorem{theorem}[stat]{Theorem}

\theoremstyle{definition}
\newtheorem{definition}[stat]{Definition}
\newtheorem{remark}[stat]{Remark}

\numberwithin{equation}{section}

\begin{document}

\title{Fractal Geometry of the Valleys of the Parabolic Anderson Equation
}
%
\author[P.\ Ghosal]{Promit Ghosal}
\address{P.\ Ghosal,
	Department of Mathematics, Massachusetts Institute of Technology (MIT),
	\newline\hphantom{\quad \ \ P. Ghosal}
	77 Massachusetts Avenue, Cambridge, MA 02139, USA
}
\email{promit@mit.edu}

\author[J.\ Yi]{Jaeyun Yi }
\address{J.\ Yi,
	Pohang University of Science and Technology (POSTECH),
	\newline\hphantom{\quad \ \ J. Yi}
	Pohang, Gyeongbuk 37673, South Korea
	}
\email{stork@postech.ac.kr}

\date{\today}

\begin{abstract} We study the macroscopic fractal properties of the \emph{deep valleys} of the solution of the $(1+1)$-dimensional parabolic Anderson equation 
\begin{equation*}
     \begin{cases} {\partial \over \partial t}u(t,x) =\frac{1}{2} {\partial^2 \over \partial x^2} u(t,x) + u(t,x)\dot{W}(t,x), \quad t>0, x\in \R,\\
   u(0,x) \equiv u_0(x), \quad x\in \R,
   \end{cases}
   \end{equation*}
where $\dot{W}$ is the time-space white noise and $0<\inf_{x\in \R} u_0(x)\leq \sup_{x\in \R} u_0(x)<\infty.$  Unlike the macroscopic multifractality of the tall peaks,  we show that valleys of the parabolic Anderson equation are macroscopically monofractal. In fact, the macroscopic Hausdorff dimension (introduced by Barlow and Taylor \cite{BT89,BT92}) of the valleys undergoes a \emph{phase transition} at a point which does not depend on the  initial data. The key tool of our proof is a lower bound to the lower tail probability of the parabolic Anderson equation. Such lower bound is obtained for the first time in this paper and will be derived by utilizing the connection between the parabolic Anderson equation and the Kardar-Parisi-Zhang equation. Our techniques of proving this lower bound can be extended to other models in the KPZ universality class including the KPZ fixed point.          
\vspace{1cm} 
 
	
\end{abstract}

\maketitle

\section{Introduction}\label{sec_intro} 
We consider the parabolic Anderson equation
\begin{equation}\label{eq:SHE}
     \begin{cases} {\partial \over \partial t}u(t,x) =\frac{1}{2} {\partial^2 \over \partial x^2} u(t,x) + u(t,x)\dot{W}(t,x), \quad t>0, x\in \R,\\
   u(0,x) \equiv u_0(x), \quad x\in \R,
   \end{cases}
   \end{equation} where $\dot{W}$ is the time-space white noise and the nonnegative initial datum $u_0\in C^{0}(\R)$ is a bounded positive initial data, i.e., 
 \begin{equation}\label{eq:condition on inital data}
   0<\inf_{x\in \R} u_0(x)\leq \sup_{x\in \R} u_0(x)<\infty.
 \end{equation}   
    The solution theory of \eqref{eq:SHE} is standard and was accomplished by Ito's calculus or martingale problem. 
    The existence and uniqueness of the solution of \eqref{eq:SHE} under those initial conditions follow from \cite[Theorem~2.2]{Bertini95} (see also \cite[Section~3.3]{Q11}). 
 Thanks to \cite{Mue91}, the solution of \eqref{eq:SHE} is strictly positive for all $t>0$ when $u_0$ is a positive initial data.
     Logarithm of the solution of \eqref{eq:SHE} formally solves the Kardar-Parisi-Zhang (KPZ) equation which is written as follows 
\begin{align}\label{eq:KPZ}
\partial_t \mathcal{H}(t,x) = \frac{1}{2} {\partial^2 \over \partial x^2} u(t,x) + \big( {\partial \over \partial x}  u(t,x)\big)^2 + \dot{W}(t,x),
\end{align}
The KPZ equation is the canonical stochastic PDE in the KPZ universality class. 
The solution theory of the KPZ equation had been approached in the recent past via different techniques, namely, the regularity structures \cite{Hai13}, paracontrolled stochastic PDE \cite{GIP15, GP17}, energy solution method \cite{gonccalves2014nonlinear}, renormalization group techniques \cite{Kup16}. The solution constructed in those works are found to be consistent with the logarithm of the solution of \eqref{eq:SHE}. The latter is the physically relevant solution of the KPZ equation and often called as the \emph{Cole-Hopf} solution.

The main objective of this paper is to study the `gaps' between the tall peaks of the parabolic Anderson equation.
The tall peaks of the solution of the parabolic Anderson equation triggers exponential growth of the moments of one point distributions. When the initial data $u_0$ is non-random and satisfies \eqref{eq:condition on inital data}, \cite{Chen15} showed that 
\begin{align}\label{eq:Intermittency}
\lim_{t\to \infty}\frac{1}{t}\log \E[u(t,x)^{k}] = \frac{k(k^2-1)}{24}, \quad \forall k\in \Z_{>0}.  
\end{align}
The above result showcases the `intermittency' (cf. \cite{CM94}) of the parabolic Anderson equation. Motivated by this result and its analogue in other stochastic PDEs with multiplicative noise, \cite{KKX17,KKX18} studied the macroscopic fractality of the spatial-temporal tall peaks of the solution for a large collection of parabolic stochastic PDEs including the parabolic Anderson equation.  They had shown that the values of the macroscopic Hausdorff dimension of the tall peaks are distinct and nontrivial when the length scale and stretch factor vary, a property which symbolizes the multifractality. This is in stark contrast with the case of Brownian motion where the tall peaks demonstrate constant Hausdorff dimension (see \cite[Theorem~1.4]{KKX17}) along different length scale.

Study on the peaks of the Parabolic Anderson models on $\Z^d$ bears many new innovations in the recent past. As we have hinted above, those works were built on the connections with the geometry of the intermittency which reveals that the total mass in parabolic Anderson model in $\Z^d$ is concentrated on smaller island of peaks which are well separated from each other. However, the mystery behind the geometry of the solution filling inner space between those islands still remains open.      
Despite of many inspiring works on the tall peaks, there was hardly any study which focuses on the so called valleys or, gaps between tall peaks. Our main goal is to showcase that the spatio-temporal valleys of the parabolic Anderson equation rather shows a different feature than the peaks and the (macroscopic) Hausdorff dimensions of the associated level sets exhibit a phase transition.


The main object of our study is the spatio-temporal level sets of the valleys shown in below 
\begin{align}
    \mathcal{V}(\gamma) : = \left\{ (t ,x) \in (e, \infty) \times \R ) \, : \, u(t,x) < e^{-\gamma t} \right\},
\end{align} 
for every $\gamma>0$. Here, $\gamma$ is called the \emph{length scale}. From \cite{CG20b}, it is known that $\log u(t,0)$ decays linearly with $t$ as $t$ grows large. In light of this fact, we may say that $\gamma$ captures the average rate of decay of the Cole-Hopf solution of the KPZ equation.
For every $\beta>0$, we define $\mathcal{S}_\beta :\R_+ \times \R \rightarrow{(1,\infty) \times \R}$ by 
\begin{align*}
    \mathcal{S}_\beta(t,x) : = \big(e^{t/\beta}, x\big) \quad \text{for all $(t,x) \in \R_+ \times \R$.}
\end{align*}
The application $\mathcal{S}_{\beta}$ on a square box produces a non-linear stretching in the time direction and the extent of stretching is determined by $\beta$ which we call as the \emph{stretch factor}.

We seek to study the fractal nature of the level sets $\mathcal{S}_{\beta}(\mathcal{V}(\gamma))$ as $t,x$ get large. The fractal nature of the peaks in parabolic Anderson equation had been quantified in \cite{KKX17,KKX18} by the Barlow-Taylor's macroscopic Hausdorff dimension. Motivated by those works, we aim to determine the macroscopic Hausdorff dimension of $\mathcal{S}_{\beta}(\mathcal{V}(\gamma))$. A precise mathematical definition of macroscopic Hausdorff dimension of any set is given in Section~\ref{subsec:dimension&localization}. For any set $A\subset \R^2$, we denote its macroscopic Hausdorff dimension by $\Dim[A]$. 

We are now ready to state the main result which will show that the macroscopic Hausdorff dimension of the valleys of the parabolic Anderson equation stays the same when we vary the stretch factor. However, it will undergo a sharp phase transition when we vary the length scale.

\begin{theorem}\label{thm:Main} Consider the solution of \eqref{eq:SHE} where $u_0$ satisfies \eqref{eq:condition on inital data}. Then, for every $\beta>0$, 
\begin{itemize}
    \item[(a)] $$\Dim[\S_\beta(\V(\gamma))] \stackrel{a.s.}{=}2 \quad \text{for } 0< \gamma < \frac{1}{24},$$
    \item[(b)] $$\Dim[\S_\beta(\V(\gamma))]\stackrel{a.s.}{=}1 \quad  \text{for }\gamma> \frac{1}{24}.$$
\end{itemize}
\end{theorem}


%
%
%
%
%
%
%
%
%
%
%

\begin{remark}
Even though Theorem~\ref{thm:Main} will be proved only for the bounded positive initial data, we believe that the same result should hold for a large class of initial data including the case when $u_0$ is a Dirac delta measure. The proof techniques will almost remain same except in few places which are pointed out Section~\ref{sec:ProofIdea}. Furthermore, it may be possible to extend Theorem~\ref{thm:lowertail} for other parabolic stochastic PDEs if few key estimates such as the ones shown in Theorem~\ref{thm:lowertail}, Proposition~\ref{prop:UpperTail} and~\ref{prop:lowertail} can be obtained in those cases.    
\end{remark}

Theorem~\ref{thm:Main} shows that for any length scale $\gamma<\frac{1}{24}$ and any stretch factor $\beta >0$, the macroscopic Hausdorff dimension of the valleys of the solution of \eqref{eq:SHE} remains constant at value $2$. The macroscopic dimension of the valleys changes to $1$ only when the stretch scale $\gamma$ is greater than $\frac{1}{24}$. This is in stark contrast with the fractal geometry of the peaks as revealed in \cite{KKX18}. On contrary to the valleys, the macroscopic Hausdorff dimension of the peaks of \eqref{eq:SHE} varies non-trivially with the stretch and length scales (see \cite[Theorem~1.1]{KKX18}). This last property marks the \emph{multifractality} of the peaks whereas Theorem~\ref{thm:Main} is a sign of the \emph{monofractality} of the valleys.

Theorem~\ref{thm:Main} in conjunction with \cite[Theorem~1.1]{KKX18} signals an apparent dichotomy between the peaks and the valleys which seems resonating with the geometry of intermittency of the parabolic Anderson model (PAM) on $\mathds{Z}^{d}$. A large number of works in the past including \cite{CM94,GKM07,KLMS09} points to fact that the most of the masses in PAM are concentrated on the very high peaks and those peaks are well separated by long stretches of trapped valleys. These distinctive features of the peaks and valleys in PAM are expected to be present in the case of parabolic Anderson equation. In fact, \cite{CCKK17} showed that if the initial data decays at infinity faster than the Gaussian kernel, then the total mass of the solution of the parabolic Anderson equation dissipates, that is, it vanishes sub-exponentially as $t\rightarrow \infty$.
 When \eqref{eq:SHE} is considered on the torus instead of $\R$, \cite{KKMS20} proved among other things that the supremum of the solution is localized in space. In combination of these works, Theorem~\ref{thm:Main} hints towards the fact that the macroscopically tall peaks of \eqref{eq:SHE} are in fact highly concentrated on the small islands which help to create many large gaps or valleys between them. 

We also like to point out that the monofractality of the valleys of \eqref{eq:SHE} does not hold in the case when the multiplicative noise $u\dot{W}$ of \eqref{eq:SHE} is replaced with the additive noise $\dot{W}$. In the latter case, the solution is a mean zero Gaussian process and \cite[Theorem~4.1]{KKX18} showed that the tall peaks of that Gaussian process are still multifractal. Due to the symmetry between the peaks and valleys for a Gaussian process, the valleys in the additive noise case are also multifractal. This attests to the fact that the monofractality of the valleys are intrinsic to the systems with multiplicative noise (see \cite{GD05,GT05,Z.et.Al00}).  

As one may see, Theorem~\ref{thm:Main} did not cover the case $\gamma = \frac{1}{24}$. Since the transition of the macroscopic dimension occurs at $\gamma = \frac{1}{24}$, we believe that one would be able to see a crossover of dimension by taking $\gamma = \frac{1}{24} + f(t)$ for some function $f:\mathds{R}_{>0}\to \mathds{R}$ such that $f(t)/t\to 0$ as $t\to 0$. More precisely, we expect the following conjecture should be true. 

\begin{conjecture}\label{conj:1}
Consider the solution of \eqref{eq:SHE} started from a bounded initial data $u_0$. Then, 
\begin{align}\label{eq:Conj}
\mathrm{Dim}_{\mathds{H}}\Big(\mathcal{S}_{\beta}\Big(\big\{(t,x)\in (e,\infty)\times \mathds{R}: u(t,t^{2/3}x)\leq e^{-\frac{t}{24} -\alpha (t\log \log t)^{1/3}}\big\}\Big)\Big)\stackrel{a.s.}{=}2-C\alpha^3
\end{align}
for all  $\alpha \in \big[0,C^{-\frac{1}{3}}\big]$ where the constant $C>0$ will depend on $\beta$ and the initial data $u_0$.
\end{conjecture}
 Conjecture~\ref{conj:1} speculates that the macroscopic dimension of $\mathcal{S}_{\beta}(\mathcal{V}(\frac{1}{24}))$ is indeed $2$ and it smoothly falls off to $1$ before $\gamma$ reaches any nonzero fixed value greater than $\frac{1}{24}$. Furthermore, the nontrivial dependence of the dimension on the stretch factor and the length scale as shown in \eqref{eq:Conj} points out to the multifractality of the valleys of \eqref{eq:SHE} at the crossover. The scaling of the spatial coordinate $x$ of $u$ by $t^{2/3}$ stems from the KPZ scaling, i.e., $1:2:3$ ratio between the scaling exponents of the fluctuation, space and time. Conjecture~\ref{conj:1} is motivated from a recent work \cite{DG21} where the authors studied the macroscopic Hausdorff dimension of the peaks of the KPZ equation at the onset of its (macroscopic) convergence towards the KPZ fixed point under the KPZ scaling. In \cite{DG21} (see the discussion after Theorem~1.3), one may find similar claim as in \eqref{eq:Conj} for the valleys. However, till now, those results can only be proved for the narrow wedge initial data of the KPZ equation which corresponds to $u_0$ being the Dirac delta measure at $0$. Proving Conjecture~\ref{conj:1} requires substantial understanding of the geometry of the KPZ equation under general initial data which we hope to explore in a future work.


One of the necessary components for studying the valleys for any stochastic process is to determine precise estimates on the probabilities of the process taking smaller values or, namely the lower tail probability. While some of the recent works indeed revealed detailed information on such estimates for the parabolic Anderson equation, the degree of the preciseness will vary with different initial data. See below for a review on those tail probability estimates. Even though those tail estimates instigate new interests in different lines of research, one of the prominent questions which will indeed play a key role in proving our results remained unanswered. This is about giving coherent lower bounds on the probabilities of the solution of \eqref{eq:SHE} taking smaller values. While such bounds are available when the initial data is a Dirac delta measure (see Theorem~1.1 of \cite{CG20b}), not much is known for any bounded initial data. Our second result of this section which we state below will partially fill this gap. 

\begin{theorem}\label{thm:lowertail}
	Suppose that $\log u_0(\cdot)$ is a bounded measurable function. For any $\gamma>\frac{1}{24}$ and $\epsilon>0$, there exists $c=c(\gamma ,\epsilon)>0$ and $t_0=t_0(\gamma, \epsilon)>0$ such that for all $t\geq t_0$,
	\begin{equation}\label{eq:LLbd}
		\P \big( u(t,x)\leq e^{-\gamma t} \big) \geq e^{-ct^{4+\epsilon}}.
	\end{equation}
\end{theorem}

\begin{remark}
We believe that the exponent of $t$ in the right hand side of \eqref{eq:LLbd} is not optimal. The lower tail large deviation principle of the KPZ equation implies that the expected exponent of $t$ is $2$. However, such result had only been proved for the case when $u_0$ is a Dirac delta measure at $0$ (see \cite{Tsai18,CC19}). \cite{CG20b} (see also Proposition~\ref{prop:LowerTail}) showed that the exponent of $t$ in the upper bound of the lower tail probability is $2$. Proving a matching exponent in the lower bound seems to be out of reach at the moment. We would also like to stress that Theorem~\ref{thm:lowertail} possibly can be extended to a large class of initial data. See Remark~\ref{rem:InitExt} for further discussion.        
\end{remark}

  Tail probabilities are ubiquitous in unlocking geometric patterns of any stochastic process. The parabolic Anderson model or, the KPZ equation are not the exceptions. Since the KPZ equation is the canonical stochastic PDE of the KPZ universality class, the techniques in handling the tail events of the KPZ equation can be replicated for other random growth models in many circumstances. To this end, the proof of Theorem~\ref{thm:lowertail} which  provides lower bound to the lower tail probability of the KPZ equation, may certainly points out to useful pathways to solve similar problems in relevant situations. Intrinsically,  Theorem \ref{thm:lowertail} ensures that the solution of \eqref{eq:SHE} cannot abruptly becomes zero. In that regards, Theorem \ref{thm:lowertail} complements the result of \cite[Theorem~1.4]{CHN16} which proves that the law of the $u(t,x)$ admits a strict positive density in $\R$. 



  Early breakthroughs in obtaining the lower tail estimates of the solution of \eqref{eq:SHE} were achieved in \cite{Mue91, MN08} using tools like large deviation bounds, comparison principles and Malliavin calculus etc. {}However, those results were proven only for bounded initial data. While \cite{Flo14} provided the first useful bound on the lower tail probability under Dirac delta initial measure,  \cite{CG20b} obtained estimates which are tight and uniform in time. The tail estimates of \cite{CG20b} were instrumental in \cite{CG20a} for further improving the tail bounds for a large class of initial data including the bounded initial data (see Proposition~\ref{prop:lowertail}). Upper tail probability of \eqref{eq:SHE} has been studied before in many places including \cite{CJK13,CD15,KKX17} in regards to its connection with the moments and intermittency property. Those bounds were recently improved in \cite{CG20a}. See Proposition~\ref{prop:UpperTail} for general initial data and Proposition~\ref{prop:uppertail} for the Dirac delta initial measure. 

Finally, we finish off this section with a few words on the novelty of the present work. One of the prime objectives of this work lies in bridging two parallel lines of research: (a) study on the fractal geometry of the PAM and related models and (b) recently discovered tools and techniques in studying long time behavior of the KPZ equation. For instance, our proof techniques which will be further elaborated in Section~\ref{sec:ProofIdea} hinge on two different approaches to study the solution of \eqref{eq:SHE}. The first one is the construction of the {\it local proxy} of the mild solution of \eqref{eq:SHE}, which is originally based on Walsh's solution theory \cite{walsh1986introduction} of stochastic PDE. The second is based on very modern tools for tackling tail probabilities of the KPZ equation. These tools are developed in last couple years by incorporating ideas from random matrix theory, geometry of random curves such as KPZ line ensemble (introduced by \cite{CH16}), interacting particle systems etc. We hope that our result will induce further interests in combining these two approaches to unravel deeper mysteries of the parabolic Anderson equation.

           
   \subsection{Proof Ideas}\label{sec:ProofIdea} 
        
   The proofs of Theorem~\ref{thm:Main} and~\ref{thm:lowertail} are mainly based on combination of probabilistic arguments and tail estimates obtained from many previous works. The core idea in the proof of Theorem~\ref{thm:Main} lies in the construction local `proxies' of the solution of \eqref{eq:SHE}. Those local objects obtained at different space-time locations will be independent of each other when the locations are far apart and in the latter case, they indeed become very close to the solution of \eqref{eq:SHE}. See Proposition~\ref{lem:localization} for more details. The construction of those local proxies are carefully carried out in \cite{KKX17,KKX18} for bounded initial data. While it could be possible to extend that construction for a large class of initial data including Dirac delta measure with some extra work, we will not pursue this direction in the present paper. Using techniques which are motivated from the works of \cite{KKX17,KKX18}, our proof will summarize the fractal information of the valleys through those local proxies. Use of the mutual independence and tail probabilities (obtained via Theorem~\ref{thm:lowertail}, Proposition~\ref{prop:lowertail},~\ref{prop:UpperTail}) of those proxies are the main takeaway tools of our analysis.   
         
 The proof of Theorem~\ref{thm:lowertail} combines several recently developed tools from \cite{CG20b,CG20a,CGH21,DG21}. The backbone of the argument is a convolution formula which provides an useful integral representation of the one point distribution of \eqref{eq:SHE} started from any reasonable data (see Proposition~\ref{prop:Convolution}). This reduces the result in Theorem~\ref{thm:lowertail} to a much simpler problem related to the solution of \eqref{eq:SHE} started from Dirac delta measure (see Theorem~\ref{prop:suplower}). We will use a fresh blend of tail estimates of the latter and monotonocity property (like FKG type inequality, see Proposition~\ref{lem:fkg}) to prove Theorem~\ref{prop:suplower} and finally, to complete proving Theorem~\ref{thm:lowertail}. The above mentioned reduction can be carried out for a large class of initial data including those which grows  linearly in the spatial variable $x$. See Remark~\ref{rem:InitExt} for more details. Our overall techniques for proving Theorem~\ref{thm:lowertail} heavily relies on the inputs from the Gibbs property of the KPZ line ensemble (via Proposition~\ref{lem:supupper} and~\ref{lem:supconti}) and monotonocity property of the KPZ equation (via Proposition~\ref{lem:fkg}). Such properties are also present in many models in the KPZ universality class including the KPZ fixed point (recently constructed in \cite{MQR16,DOV18}), asymmetric simple exclusion process (ASEP), last passage percolation, stochastic six vertex model etc. We hope that it could be possible to prove lower bound to the lower tail probability of those models using the techniques of the present paper.             

\subsection*{Outline}
The rest of the paper is organized as follows. Section~\ref{sec_prelim} reviews some preliminary facts and estimates on the tail probabilities of \eqref{eq:SHE} and furthermore, recalls the definition and some useful results on the macroscopic Hausdorff dimension. Theorem~\ref{thm:Main} will be proved in Section~\ref{sec_proofthm}. Finally, Section~\ref{sec:Tail} will contain the proof of Theorem~\ref{thm:lowertail}.             

\section{Preliminaries} \label{sec_prelim}

In this section, we introduce few notations and recall some known facts which are important for our analysis. 

\subsection{Convolution Formula \& Tail Estimates}\label{subsec:convolution&tail}
Recall that the logarithm of the solution of \eqref{eq:SHE} is the Cole-Hopf solution of the KPZ equation. When the initial data $u_0$ is a Dirac delta measure at $0$, we denote the Cole-Hopf solution by $\mathcal{H}^{\mathbf{nw}}$ where `$\mathbf{nw}$' stands for the narrow wedge initial data for the KPZ equation. 

 The solution of \eqref{eq:SHE} started from the Dirac delta measure is called the fundamental solution. The one dimensional distributions of $u(t,x)$ started from any initial data can be descried in terms of the fundamental solution via the \emph{convolution formula}. The precise statement is stated as follows. 
 
\begin{proposition}[Convolution Formula, Lemma 1.18 of \cite{CH16}]\label{prop:Convolution}
	For any measurable function $u_0: \R\to \R_{>0}$ and for a fixed $(t,x)\in \R_{\geq 0} \times \R$, the unique solution of \eqref{eq:SHE} started from $u_{0}$ satisfies
	\begin{align}\label{eq:Convolution}
		\log u(t,x)\deq  \log \left( \int_{-\infty}^\infty u_0(x-y)e^{\mathcal{H}^{\mathbf{nw}}(t,y)}dy\right)
	\end{align} 
\end{proposition}  

In some of the recent works \cite{CG20a,GL20}, the convolution formula has been utilized to derive useful bounds on the probability of $u(t,x)$ getting too large or, too small. In the next two propositions, we summarize those findings. For the sake of completeness, we will provide a short proof of these propositions in Section~\ref{sec:UpperTail} and~\ref{sec:LowerTail} respectively.  

\begin{proposition}[Upper Tail Estimates]\label{prop:UpperTail} 
Fix $\gamma <\frac{1}{24}$. There exist $c_1>c_2>0$ and $t_0>0$ which depend on $\gamma$ and initial data $u_0$ such that for all $t>t_0$ and $x\in \R$,  
\begin{equation}\label{eq:lya}
   e^{ -c_1(\frac{1}{24} -\gamma)^{3/2} t}\leq  \P(u(t,x)>e^{-\gamma t})\leq e^{ -c_2(\frac{1}{24} -\gamma)^{3/2} t}. 
\end{equation} 
\end{proposition}

\begin{proposition}[Lower Tail Estimates]\label{prop:lowertail} Fix any small $\delta>0$ and $\gamma>\frac{1}{24}$. There exist  $c= c(\gamma,\delta,u_0)>0$ and $t_0(\gamma,\delta,u_0)>0$ such that for any $t\geq t_0$ and $x\in \R$, 
\begin{equation}\label{eq:lowertail}
    \P \big( u(t,x) < e^{-\gamma t} \big) \leq e^{ -c(\gamma - \frac{1}{24})^{3-3\delta/2}t^{2-\delta}}. 
\end{equation} 
         \end{proposition}

\begin{remark}
Notice that \eqref{eq:lowertail} found an upper bound to the lower tail probability. As we have mentioned before, the exponent of $t$ in the upper bound of \eqref{eq:lowertail} and the lower bound of \eqref{eq:LLbd} are not close.   
\end{remark}

\subsection{Macroscopic dimension \& Localization}\label{subsec:dimension&localization} 

In this section, we recall the definition of Barlow-Taylor's macroscopic Hausdorff dimension in $\R_{+}\times \R$ and and some of its properties. Localization technique in parabolic Anderson model, pioneered by \cite{KKX17} is an important tool for studying the macroscopic fractality. We review this tool in Proposition~\ref{lem:localization}.     

\begin{definition}[Barlow-Taylor's macroscopic Hausdorff dimension \cite{BT89,BT92}]\label{def:HausdorffDim}
Let $\mathcal{B}_r$ be the collection of all sets of the form 
\begin{align}
Q(x,r) := [x_1,x_1+r)\times [x_2,x_2+r)
\end{align}
for $x:= (x_1,x_2)\in \R_{+}\times \R$ and $r$ in $(0, \infty)$.  Define $\mathcal{B} := \bigcup_{r\geq 1} \mathcal{B}_r$. Let $\mathds{V}_n:= [0,e^{n})\times [-e^n,e^n)$ for $n\in \Z_{\geq 0}$ and $\mathds{S}_n:= \mathds{V}_n\backslash \mathds{V}_{n-1}$ for $n \in \Z_{>0}$. For any $E\subset \R^d$ and $\rho>0$, the $\rho$-dimensional macroscopic Hausdorff content of the set $E$ will be denoted as $\nu_{n,\rho}(E)$ and defined as 
\begin{align}
\nu_{n,\rho}(E) := \inf \sum_{i=1}^{m}\Big(\frac{\mathrm{side}(Q_i)}{e^n}\Big)^{\rho}
\end{align}
where the infimum is taken over all $Q_1, \ldots , Q_m \in \mathcal{B}$ such that $\cup_{i=1}^{m}Q_i$ covers the set $E\cap \mathds{S}_n$. The Barlow-Taylor's \emph{macroscopic Hausdorff dimension} of any set $E$ which we will denote as $\Dim(E)$ is defined as 
\begin{align}
\Dim(E) := \inf\Big\{\rho>0: \sum_{n=1}^{\infty}\nu_{n,\rho}(E)<\infty\Big\}.
\end{align}
\end{definition}

The following proposition is useful in getting lower bound to the macroscopic Hausdorff dimension. 

\begin{proposition}[Theorem~4 of \cite{BT92}]\label{prop:DensityTheorem}
Fix $\gamma \in (0,2)$. For any set $E$ and $n\in \Z$, let us define  
\begin{align}\label{eq:mu_n}
\mu_n(E) = \sum_{\substack{s\in \Z\\e^{n}< s\leq e^{n+1}}} \sum_{0\leq j < e^{n(1-\gamma)}} \mathbf{1}((s,j)\in E). 
\end{align} 
Then, there exists a constant $C>0$ such that $\nu_{n,2-\gamma}(E)\geq C e^{-(2-\gamma)n}\mu_n(E)$. 
\end{proposition}

The next result which is taken from Corollary~1.4 of \cite{KKX18} constructs a two dimensional set whose macroscopic Hausdorff dimension is $1$. We use this set to show the phase transition of the macroscopic Hausdorff dimension in Theorem~\ref{thm:Main}.

\begin{proposition}\label{prop:SpSetDim}
Fix an arbitrary constant $q>1$ and for any $n\geq 1$, define 
\begin{align}\label{eq:Xi}
\Dim(\Xi_q) = 1, \quad \Xi_q:= \Big\{(x,y)\in (0,\infty)^2: y\geq x^q\Big\}.  
\end{align} 
\end{proposition} 
Corollary~1.4 of \cite{KKX18} also showed that the dimension of $\Xi_q$ is $2$ for any $q\in (0,1]$.   
  
The next result describes an important tool (originated from \cite{CJK13}) for studying the fractal geometry pioneered by \cite{KKX17, KKX18}. In brief, it describes how to construct `local proxy' for the the solution of \eqref{eq:SHE} in far away spots so that the proxies will be independent of each other. In \cite{KKX17, KKX18}, these proxies were made to summarize the fractal information hidden in the peaks of the parabolic Anderson equation. As we will show in our proof, these proxies are also useful in decoding fractality of the valleys.     
    
\begin{proposition}[`Local Proxy' in parabolic Anderson equation, Theorem~3.9 of \cite{KKX18}] \label{lem:localization}
Fix $k\geq 2$. There exists a finite constant $c_0>0$ independent of $k$ and $t$ such that for any finite set of nonrandom points $x_1,...,x_m \in \R$ and $t\geq1$ that satisfy 
\begin{equation}\label{localization1}
    \min_{1\leq i\neq j \leq m } |x_i-x_j| > c_0t^2 k^3,
\end{equation}
there exists a set of independent random variables $Y^{(k)}_1, Y^{(k)}_2, \ldots , Y^{(k)}_m $ with all positive moments finite, satisfying  
\begin{equation}\label{eq:error:u-uk}
   \sup_{1\leq j\leq m} \E\Big[ \big|u(t,x_j)-Y_j^{(k)}\big|^k  \Big] \leq C^k e^{- k^3 t},
\end{equation}
where $C$ is a positive constant independent $k$, $t$ and $\{x_i:1\leq i\leq m\}$. 
\end{proposition}

%



\section{Fractality of Valleys} \label{sec_proofthm}
\subsection{Proof of part (a) of Theorem \ref{thm:Main}}\label{subsec_part1}
Since $\mathcal{S}_{\beta}(\mathcal{V}(\gamma))$ is subset of a $\R_{+}\times \R$,  
$\Dim(\mathcal{S}_{\beta}(\mathcal{V}(\gamma)))\leq 2$. It suffices to show that $\Dim(\mathcal{S}_{\beta}(\mathcal{V}(\gamma)))\geq 2$. The definition of the Barlow-Taylor's Hausdorff dimension implies that if $\sum_{n=1}^{\infty}\nu_{n,2-\gamma}(\mathcal{S}_{\beta}(\mathcal{V}(\gamma)))=\infty$ holds with probability $1$ for all $\gamma$, then, $\Dim(\mathcal{S}_{\beta}(\mathcal{V}(\gamma)))$ is almost surely greater than or equal to $2$.     In what follows, we will show that this indeed holds, i.e., 
\begin{align}\label{eq:nu_sum}
\sum_{n=1}^{\infty}\nu_{n,2-\gamma}(E) \stackrel{a.s.}{=} \infty, \quad \forall \gamma \in (0,2).  
\end{align}
Our proof will require Proposition~\ref{prop:SpSetDim} and~\ref{prop:UpperTail}. The  first proposition which is taken from \cite{BT92} bounds $\nu_{n,2-\gamma}$ from below by a discrete measure (see \eqref{eq:mu_n}) on the set $(e^{n}, e^{n+1}]\times (0, e^{n(1-\gamma)}$ for any $\gamma \in (0,2)$. The main goal of this section will be to obtain suitable lower bound to that discrete measure. One of the key tools that we use is the upper tail probability of the parabolic Anderson equation. Tight bounds to this tail probability was derived in \cite{CG20a} for a large class of initial data and we have compiled their result in Proposition~\ref{prop:UpperTail}.   
We now proceed to complete the proof of part $(a)$ of Theorem~\ref{thm:Main}.


\noindent \emph{Proof of Theorem~\ref{thm:Main}(a):} Consider a set of points $x_1,\ldots ,x_m\in \R$ satisfying \eqref{localization1}. We seek to show that $\mu_n(S_{\beta}(\mathcal{V}(\gamma)))$ is bounded below by $e^{(2-\gamma)n}$ almost surely for all large $n$. This will show \eqref{eq:nu_sum} and hence, will complete proof. We divide the rest of the proof into two steps. In \emph{Step I}, we will show a upper bound to the tail probability of $\max_{1\leq j\leq m} u(t,x_j)$. \emph{Step II} will use this upper bound to show that $\mu_n(S_{\beta}(\mathcal{V}(\gamma)))$ is almost surely close to $e^{(2-\gamma)n}$ for all large $n$.
\smallskip 

\noindent \emph{Step I:} Fix $t_0>0$ and $\gamma \in (0,\frac{1}{24})$. Here, we claim and prove that there exists $\mu=\mu(t_0,\gamma)>0$, $C= C(t_0,\gamma)>0$ and $k$ large such that for all $t>t_0$,
    \begin{align}\label{eq_mintail}
        \P \Big( \min_{1\leq j\leq m} u(t,x_j) > e^{-\gamma t} \Big)\leq  e^{-m\mu t +m\log 2  } + C^kme^{(-k^3 + \gamma k)t}. 
\end{align}

By Proposition~\ref{prop:UpperTail}, for any $\gamma<\frac{1}{24}$ and $t_0>0$, there exists $\mu:=\mu(t_0,\gamma)>0$ such that for all $t>t_0$ and $x\in \R$,
\begin{equation} 
    \P\big( u(t,x) > e^{-\gamma t} \big) \leq e^{-\mu t}.
\end{equation}  

Recall the local proxy $\{Y^{(k)}_j\}_j $ defined in Proposition~\ref{lem:localization}. By the union bound, we get 
\begin{align*}
    \P \big( \min_{1\leq j\leq m} u(t,x_j) > e^{-\gamma t} \big) \leq &\P\big(  \min_{1\leq j\leq m} Y^{(k)}_j > \frac{e^{-\gamma t}}{2}\big) +\P \big(\max_{1\leq j\leq m} |u(t,x_j)-Y^{(k)}_j | >  \frac{e^{-\gamma t}}{2}\big).
\end{align*}
In what follows, we bound the two terms on the right side of the above display. By Lemma~\ref{localization1}, for any $k>2$, $\{Y^{(k)}_j\}_{j}$ is a set of independent random variables which imply
\begin{align}
    \P\big(  \min_{1\leq j\leq m} Y^{(k)}_j > \frac{1}{2}e^{-\gamma t}\big) &= \P\big(Y^{(k)}_j > \frac{1}{2}e^{-\gamma t}\big)^m \nonumber\\
    &\leq \Big( \P \big( u(t,x) > \frac{1}{4}e^{-\gamma t} \big) +\P \big( |u(t,x_j) - Y^{(k)}_j | > \frac{1}{4}e^{-\gamma t} \big) \Big)^m\nonumber\\
    &\leq \Big( e^{-\mu t } + C^ke^{(-k^3 + k\gamma )t} \Big)^m \leq 2^m e^{-m\mu t }. \label{eq:1stTermBd}
\end{align}
The last inequality holds for $k$ large since $e^{-\mu t}> C^k e^{-(k^3-k\gamma)t}$ for large $k$, uniformly in $t$.  
Moreover, the union bound and Lemma~\ref{localization1} yields
\begin{align}
    \P \big(\max_{1\leq j\leq m} |u(t,x_j)-Y^{(k)}_j | >  e^{-\gamma t}\big) &\leq \sum_{j=1}^m \P \big( |u(t,x_j) - Y^{(k)}_j| > e^{-\gamma t}\big)\nonumber\\ &\leq mC^k e^{(-k^3 +\gamma k) t}.\label{eq:2ndTermBd}
    \end{align} 
Combining \eqref{eq:1stTermBd} and \eqref{eq:2ndTermBd} yields \eqref{eq_mintail}. 
\smallskip



\noindent \emph{Step II:}
Define the random set 
\begin{align*}
    \mathfrak{S}_{\beta,\gamma}= \left\{(t,x) \in (1,\infty)\times \R \,:\, u(\beta \log t , x) < t^{-\beta \gamma} \right\},\quad  \hat{\mathfrak{S}}_{\beta,\gamma} := \mathfrak{S}_{\beta,\gamma} \cap \bigcup_{n=0}^\infty \left( e^n,e^{n+1}\right]^2.
\end{align*}
 Notice that $\mathfrak{S}_{\beta, \gamma}(\gamma)$ is same as $ \{(t,x): (\beta \log t,x)\in \mathcal{S}_{\beta}(\mathcal{V}(\gamma))\}\cap [1,\infty)\times \R$. 
 Fix $0<\delta<\epsilon<1$. For any $n\in \N$  and   $j\in [0,e^{n(1-\delta)})\cap \Z$, define 
\begin{align*}
    a_{j,n}:= e^n + je^{n\delta}, \qquad  \I_{j,n}(\delta) := (a_{j,n}, a_{j+1,n}].
\end{align*}  
 Our goal is now to obtain an upper bound for the supremum of the probabilities of the events $\{\inf_{x\in \mathcal{I}_{n,j}} u(\beta \log t,x)>t^{-\beta \gamma}\}$ as $j$ varies over the set $[0,e^{n(1-\delta)})\cap \Z$. 
 For this, we seek to apply \eqref{eq_mintail}.   Pick an integer $m$ from the set  $\Big[\frac{1}{2}e^{n(\delta-\epsilon)}, 2e^{n(\delta-\epsilon)}\Big]\cap \Z$. Fix a set of integers $x_1,...,x_m \in \I_{j,n}(\delta)$ such that for all $n>1$,
    $  \min_{1\leq i\neq j \leq m} |x_i-x_j| \geq e^{n \epsilon}$. 
As a consequence, the following inequality holds 
\begin{equation}
    \min_{1\leq i\neq j \leq m} |x_i-x_j| \geq e^{n \epsilon} \geq c_0 k^3 (\beta \log t )^2
\end{equation} for $t\in \left( e^n , e^{n+1}\right]$ when $n$ is sufficiently large. Since $\inf_{s\in \mathcal{I}_{n,j}} u(\beta \log t ,x)$ is less than the infimum value $ u(\beta \log t ,x_i)$ for $1\leq i\leq m$, we may write  
\begin{align}
     \P \Big( \inf_{s\in \mathcal{I}_{n,j}} u(\beta \log t ,x) >t^{-\gamma \beta}\Big) & \leq \P \Big( \inf_{1\leq i\leq m} u(\beta \log t ,x_i) >t^{-\gamma \beta}\Big) \nonumber \\ &\leq \exp\left( -\frac{1}{2}e^{n(\delta-\epsilon)} \left(\mu \beta n - \log 2  \right) \right)+ C^k\exp\left( (\delta-\epsilon-(k^3-\gamma k )\beta )n \right),\label{eq_maxesti}
\end{align} for all $t\in (e^n,e^{n+1}]$ with all sufficiently large $n$. The last inequality follows by applying \eqref{eq_mintail} and recalling that $2^{-1}e^{n(1-\epsilon)}\leq m \leq 2e^{n(1-\epsilon)}$. Note that right hand side of the above inequality does not depend on $j$. 
Therefore, we can deduce the following inequality
\begin{align*}  
    \P &\Big( \hat{\mathfrak{S}}_{\beta,\gamma}  \cap (\{t\} \times \I_{j,n}(\delta)) =\varnothing \text{ for some } j \in [0, e^{n(1-\delta)}) \cap \Z \text{ and }t\in (e^n,e^{n+1}]\cap \Z \Big)\\
    &\leq  \P \Big(   \Big\{\max_{\substack{t\in \Z \\ t\in(e^n,e^{n+1}]}}\max_{\substack{j\in \Z \\ j\in(0,e^{n(1-\delta)}] } } \inf_{x \in \I_{j,n}(\gamma) } t^{\gamma \beta }u(\beta \log t ,x )\Big\} >1 \Big)\\
    &\leq \sum_{t\in \Z\cap (e^n,e^{n+1}]} \sum_{j \in \Z\cap (e^n,e^{n+1}]} \P \Big( \inf_{s\in \mathcal{I}_{n,j}} u(\beta \log t ,x) >t^{-\gamma \beta}\Big) 
     \\&\leq e^{n(2-\delta)+1} \left[\exp\left( -\frac{1}{2}e^{n(\delta-\epsilon)} \left(\mu \beta n - \log 2  \right) \right)+ C^k\exp\left( (\delta-\epsilon-(k^3-\gamma k )\beta )n \right)\right].
\end{align*} 
The last inequality follows from the direct application of \eqref{eq_maxesti}. For large $k$, the right hand side of the last inequality is summable in $n$. Thus, by the Borel-Cantelli lemma, the following holds
\begin{equation}
     \hat{\mathfrak{S}}_{\beta,\gamma} \cap (\{t\} \times \I_{j,n}(\delta)) \neq\varnothing \text{ for all } j \in [0, e^{n(1-\delta)}) \cap \Z \text{ and }t\in (e^n,e^{n+1}]\cap \Z,
\end{equation} almost surely for all sufficiently large $n$. As a result, $\mu_{2-\delta}(\mathcal{S}_{\beta}(\mathcal{V}(\gamma)))$ eventually exceeds $Ce^{n(2-\delta)}$ as $n$ increases implying $\sum_{n}\mu_{n,2-\delta}(\mathcal{S}_{\beta}(\mathcal{V}(\gamma)))= \infty$ holds almost surely. Combining this with Proposition~\ref{prop:DensityTheorem} yields $\Dim(\mathcal{S}_{\beta}(\mathcal{V}(\gamma)))\geq 2-\delta$ almost surely. Letting $\delta $ to $0$ competes the proof.

\qed
\subsection{Proof of part (b) of Theorem \ref{thm:Main}}

We complete the proof in two steps. The first step is to show that $\Dim[\S_\beta(\V(\gamma))]$ is bounded above by  $1$ and the second step is to show that $\Dim[\S_\beta(\V(\gamma))]$ is bounded below by $1$. These two steps will be contained in the following two propositions.
These propositions together complete the proof of part $(b)$. 
\begin{proposition}\label{prop:dimupper} For all $\gamma >1/24$ and $\beta >0$, we have
\begin{equation}\label{eq:DimUpper}
       \Dim[\S_\beta(\V(\gamma))]=  \Dim \Big(\Big\{ (t,x) \in (1,\infty) \times \R : u(\beta\log t, x) < t^{-\beta \gamma} \Big\}\Big)\leq 1 \quad \text{a.s.}
    \end{equation}    
\end{proposition}

\begin{proposition}\label{prop:dimlower}
	For all $\gamma > 1/24$ and $\beta>0$, there exists $t_0>0$ such that for all $t\geq t_0$, 
	\begin{equation}\label{eq:DimLow}
		\Dim[\S_\beta(\V(\gamma))]\geq 	\Dim\left(\left\{ x\in \R : u(\beta \log t , x ) < t^{-\beta \gamma } \right\}\right) \geq 1 \quad \text{a.s.}
	\end{equation}
\end{proposition}

One of the key inputs for the proof of these propositions is the upper bound on the probability of $u(t,x)$ being small in a bounded domain. This result will be obtained in Proposition~\ref{prop:inf} using Proposition \ref{prop:lowertail}.

\begin{proposition}\label{prop:inf} Let $\gamma >1/24$ and $\epsilon\in(0,1)$. There exists a positive constant $C=C( \gamma,\epsilon)$ such that for all large $a\geq 1$, $b>0$, and integers $l_1,l_2\geq 1 $, 
 \begin{align}\label{eq:inf}
    \P \big(\text{For some }(t,x) \in (a,a+l_1] \times (b,b+l_2], u(t,x) <e^{-\gamma t} \big)
    \leq Cl_1l_2e^{-C a^{2-\epsilon}}.
\end{align}
\end{proposition}

\begin{remark}
  \cite[Theorem~1.4]{CK17} established an upper bound on the lower tail probability of the infimum of $u(t,x)$ being smaller than fixed $\epsilon>0$. Proposition~\ref{prop:inf} shows that a tighter upper bound \eqref{eq:inf} holds when $t$ gets larger.  
\end{remark}

We first prove Proposition~\ref{prop:dimupper} and~\ref{prop:dimlower} assuming Proposition~\ref{prop:inf} in two ensuing subsections and then, Proposition~\ref{prop:inf} will be finally proved in Section~\ref{subsec:proof oof prop:inf}.

\subsubsection{Proof of Proposition~\ref{prop:dimupper}} \label{subsec:proof of prop:dimupper}

The equality in \eqref{eq:DimLow} between the macroscopic Hausdorff dimension of $\S_\beta(\V(\gamma))$ and the set $\{(t,x) \in (1,\infty) \times \R : u(\beta\log t, x) < t^{-\beta \gamma}\}$ is straightforward. 

We now prove that the dimension of the latter set is at most $1$. Observe that for all $a\geq 1 $, $\beta \log(a+1) - \beta\log a = \beta \log\big( 1+\frac{1}{a}\big)\leq \beta \log 2 $. Fix an arbitrary constant $q >1$. Proposition~\ref{prop:inf} implies that for any $\delta>0$, $\gamma >1/24$, $\beta,b>0$, there exists $C=C(\delta, \gamma)>0$ such that uniformly for all $a\in(e^{n/q},e^{n+1}]$ with all large $n$,   
    \begin{align}
        \P &\Big( \text{For some }t\in (a,a+1]:\inf_{x\in (b,b+1]} u(\beta \log t ,x) < t^{-\beta \gamma})\Big)\leq C\exp\Big(-C ( \beta n/q)^{2-\delta} \Big)\label{eq:infupperq}
    \end{align}     
    

  Below we introduce few notations which will be used throughout the rest of the proof. Define $\mathcal{I}_{n} := (e^{n},e^{n+1}]$, $\mathcal{I}^{(q)}_n:=(e^{n/q}, e^{n+1}]$, $\mathcal{J}^{(q)}_n:= (0,e^{n/q}]$ and
  \begin{equation}\label{eq:GbDef} 
 \begin{aligned}
        \G_{\beta, \gamma} : = \Big\{ (t,x) \in (1,\infty) \times (0,\infty) &\,:\, u(\beta \log t ,x) > t^{-\beta \gamma}\Big\},\\
  \L_n : = \L_n(\beta, \gamma) : = \G_{\beta, \gamma } \cap (\mathcal{I}^{(q)}_n\times \mathcal{I}^{(q)}_n), \quad & \quad   \L'_n : = \big(\mathcal{I}_n \times \mathcal{J}^{(q)}_n\big) \cup \big(\mathcal{J}^{(q)}_n \times \mathcal{I}_n\big).     
\end{aligned} 
\end{equation}

We now claim and prove that 
\begin{align}\label{eq:GDim}
         \Dim (\G_{\beta ,\gamma}) \leq 1 \quad \text{a.s.}
     \end{align}
Before proceeding to the proof of \eqref{eq:GDim}, let us explain how \eqref{eq:GDim} completes the proof of the inequality in \eqref{eq:DimUpper}. As one will see below, similar argument as in the proof of \eqref{eq:GDim} will imply  
\begin{align}
 \Dim \Big( \big\{ (t,x) \in (1,\infty)\times (-\infty ,0): u(\beta\log t ,x) < t^{-\beta\gamma}\big\} \Big) \leq 1.
\end{align} Moreover, we have $\Dim( (1,\infty) \times \{0\})=1$ (see \cite[Section~4.1]{BT92}). Combining these with \eqref{eq:GDim} and recalling that $\mathrm{Dim}_{\mathds{H}}(A\cup B) = \max\{\mathrm{Dim}_{\mathds{H}}(A), \mathrm{Dim}_{\mathds{H}}(B)\}$ for any two sets $A,B\subset \mathbb{R}^2$ show the upper bound in the macroscopic Hausdorff dimension from \eqref{eq:DimUpper}.

Now we proceed to prove \eqref{eq:GDim}. Recall $\mathds{S}_n$ from Definition~\ref{def:HausdorffDim}. From \eqref{eq:GbDef}, it follows
    \begin{equation} \label{eq:decomp}   
                \G_{\beta, \gamma } \cap \mathds{S}_{n+1} \subseteq \L_n \cup \L'_n.
                    \end{equation}       
     Using Proposition~\ref{prop:SpSetDim}, we can deduce that $\Dim \left( \cup_{n=1}^\infty \L'_n\right) \leq 1 $. From this fact and \eqref{eq:decomp}, we have 
     \begin{equation}\label{eq:dimLn}
         \Dim( \G_{\beta , \gamma}) \leq \max \Big\{ 1, \Dim \big(\bigcup_{n=1}^\infty \L_n \big)\Big\} \quad \text{a.s.}
     \end{equation} In the view of \eqref{eq:dimLn}, we now bound the dimension of $\cup_{n=1}^\infty \L_n$. Note that we can cover $\mathcal{I}^{(q)}_n\times \mathcal{I}^{(q)}_n$ with $O(e^{2n})$ many squares of the form $(a,a+1]\times (b,b+1]$ for all integers $n\geq 1 $. Let us denote the set of all these squares needed to cover $\mathcal{I}^{(q)}_n\times \mathcal{I}^{(q)}_n$ by $\Theta_{n,q}$. Out of those $O(e^{2n})$ squares in $\Theta_{n,q}$, $\L_n$ can be covered with those satisfying 
     \begin{equation}\label{condi:squares}
         \inf_{x\in (b,b+1]} u(\beta \log t , x ) < t^{-\beta \gamma }, \text{for some } t\in(a,a+1],. 
     \end{equation}      
     Thus, for all large integers $n\geq 1 $ and for all real $\rho>0$, 
     \begin{align*}
         \E \big[ \nu_{n,\rho}(\L_n)\big]&\leq e^{-n\rho}\sum_{a,a+1]\times (b,b+1]\in \Theta_{n,q}}\P \big(\inf_{x\in (b,b+1]} u(\beta \log t , x ) < t^{-\beta \gamma }, \text{for some } t\in(a,a+1]\big) \\
         &\leq C \exp \big(  (2-\rho)n -C \big(\beta n/q\big)^{2-\delta }  \big). 
     \end{align*} 
where the last inequality follows from \eqref{eq:infupperq}.      
     Here, the positive constant $C$ is independent of $(\rho, n,q)$. Note that the right hand side of the above display geometrically decays with $n$. Summing both sides of the above display over $n$, we deduce that $\sum_{n=1}^\infty \nu_\rho^n(\L_n) <\infty$ a.s., for every $\rho >0$. From the definition of the macroscopic Hausdorff dimension, this implies that $ \Dim \left(\bigcup_{n=1}^\infty \L_n \right) =0$ a.s. Together with \eqref{eq:dimLn}, the above result completes the proof of \eqref{eq:GDim}. 

\subsubsection{Proof of Proposition~\ref{prop:dimlower}} \label{subsec:proof of prop:dimlower}

	Note that 
	\begin{equation*}
		\left\{ (s,x) \in (1,\infty) \times \R : u(\beta\log s, x) < s^{-\beta \gamma} \right\} \supseteq \{t\}\times \left\{ x\in \R : u(\beta \log t , x ) < t^{-\beta \gamma } \right\},
	\end{equation*} for all $t \geq 1 $. From the above display, the first inequality of \eqref{eq:DimLow} follows via the monotonocity of the Hausdorff dimension. Now we proceed to prove the second inequality. As in Section \ref{subsec_part1}, we define the following random sets for any $t>1$,
	\begin{align*}
 \mathfrak{S}^{(t)}_{\beta,\gamma}:= \left\{(t,x) \in (1,\infty)\times \R \,:\, u(\beta \log t , x) < t^{-\beta \gamma} \right\},\quad  \hat{\mathfrak{S}}^{(t)}_{\beta,\gamma} := \mathfrak{S}^{(t)}_{\beta,\gamma} \cap \bigcup_{n=0}^\infty \left( e^n,e^{n+1}\right].	
\end{align*}	

Let us choose and fix arbitrary reals $\delta$ and $\epsilon$ satisfying $0<\epsilon < \delta < 1 $. For any  $n\in \Z_{\geq 0}$ and $j\in [0,e^{n(1-\delta)})\cap \Z_{\geq 0}$, recall the definitions  $a_{j,n}:= e^n + je^{n\delta}$ and $\I_{j,n}(\delta) := (a_{j,n}, a_{j+1,n} ].$ 
In what follows, we claim and prove that for all sufficiently large $t_0=t_0(\gamma)>1$ , there exist $c_1 = c_1(t_0)>0, c_2 = c_2(t_0)>0$ and $C = C(t_0)>0$ such that for all $t\geq t_0$,  
\begin{align}\label{eq:NDelta} 
    \max_{j \in [0,e^{n(1-\delta)}) \cap \Z} &\P \Big( \inf_{x \in \I_{j,n}(\delta) } u(\beta \log t ,x) >t^{-\gamma \beta}\Big) \\&\leq \exp\left(-c_1 e^{c_2t^{4}} e^{n(\delta -\epsilon )}\right) + C^n \exp\left( n(\delta-\epsilon) -n^3t +n(\log 2 +\gamma t) \right).
\end{align}
 
Let us fix $x_1,...,x_m \in \I_{j,n}(\delta)$ such that for all $n>1$, the following conditions are satisfied:
 \begin{align}  
 \min_{1\leq i\neq j \leq m} |x_i-x_j| \geq e^{n \epsilon}, \qquad \frac{1}{2}e^{n(\delta-\epsilon)}\leq m\leq 2e^{n(\delta-\epsilon)}.
 \end{align}
 We first bound the probability of $\min_{1\leq j\leq m} u(t,x_j)$ exceeding the value $e^{-\gamma t}$. We will use this to prove \eqref{eq:NDelta}.

 Recall $Y^{(k)}_j$ with $k=n$ from Proposition~\ref{lem:localization}. By the union bound, we get 
\begin{equation}\label{eq:minesti}
	\begin{aligned}
    \P \Big( \min_{1\leq j\leq m} u(t,x_j) > e^{-\gamma t} \Big) \leq &\P\Big(  \min_{1\leq j\leq m} Y^{(n)}_j > \frac{1}{2}e^{-\gamma t}\Big) \\ &+ \P \Big(\max_{1\leq j\leq m} |u(t,x_j)-Y^{(n)}_j | >  \frac{1}{2}e^{-\gamma t}\Big).
\end{aligned}
\end{equation}
In what follows, we bound two terms of the right hand side of the above display. Recall that $\{Y^{(n)}_j\}_{j=1}^{m}$ are independent which implies   
\begin{align}\label{eq:EqString}
    \P\Big(\min_{1\leq j\leq m} Y^{(n)}_j > \frac{1}{2}e^{-\gamma t}\Big) &= \P\Big(   Y^{(n)}_j > \frac{1}{2}e^{-\gamma t}\Big)^m \nonumber\\
    &\leq \Big( \P\big(u(t,x) > \frac{1}{4}e^{-\gamma t}\big)+\P \big( |u(t,x_j) - Y^{(n)}_j | >\frac{1}{4} e^{-\gamma t} \big) \Big)^m \nonumber
    \\
    &\leq \Big(1- \P \big( u(t,x)<\frac{1}{4}e^{-\gamma t} \big)+\P \big( |u(t,x_j) - Y^{(n)}_j | > \frac{1}{4}e^{-\gamma t} \big) \Big)^m.
\end{align}

By Theorem \ref{thm:lowertail}, for any $t_0>\frac{\log4}{\gamma}$ and $\epsilon'>0$, there exists $c_1=c_1(t_0,\epsilon')>0$ and $c_2= c_2(t_0, \epsilon')>0$ such that for all $t\geq t_0$ 
\begin{equation}\label{eq:A}
	\P \big( u(t,x) < \frac{1}{4}e^{-\gamma t } \big) \geq \P \big( u(t,x) < e^{-2\gamma t } \big)\geq  \exp(-ct^{4+\epsilon'}).
\end{equation}  We fix $t\geq t_0$ and use Chebyshev's inequality in conjunction with Proposition \ref{lem:localization} to conclude  
\begin{equation}\label{eq:B}
	\P \Big( |u(t,x_j) - Y^{(k)}_j | > \frac{1}{4}e^{-\gamma t} \Big) \leq C^n e^{-n^3 t +n(\log 4+\gamma t) }\leq \frac{1}{2}\exp(-ct^{4+\epsilon'}),
\end{equation} 
where the last inequality follows for all large $n$. 
 Combining \eqref{eq:A}, \eqref{eq:B} and substituting those into the last line of \eqref{eq:EqString} shows 
\begin{equation*}
	\P\Big(  \min_{1\leq j\leq m} Y^{(k)}_j > \frac{1}{2}e^{-\gamma t}\Big) \leq \big( 1-\frac{1}{2}\exp(-ct^{4+\epsilon'})\big)^m \leq \exp\left( -\frac{m}{2}e^{-ct^{4+\epsilon'}}   \right),
\end{equation*}
where the last inequality holds since $1-x \leq e^{-x}$ for all $x>0$. This provides an upper bound to the first term in the right hand side of \eqref{eq:minesti}. The second term can be bounded above by using the first inequality of \eqref{eq:B} and union bound. As a result, we obtain the following  
\begin{equation}\label{eq:mincontrol}
	\text{r.h.s. of \eqref{eq:minesti}} \leq \exp\left( -\frac{m}{2}e^{-ct^{4+\epsilon'}}   \right) +  mC^n e^{-n^3 t +n(\log 2+\gamma t) }.
\end{equation}  Hence, 
\begin{align*}
    \P &\Big( \mathfrak{S}^{(t)}_{\beta, \gamma} \cap   \I_{j,n}(\delta) =\varnothing \text{ for some } j \in [0, e^{n(1-\delta)}) \cap \Z \Big)\\
    \leq & \P \Big(  \max_{j\in(0,e^{n(1-\delta)}]\cap \Z} \inf_{x \in \I_{j,n}(\gamma) } t^{\gamma \beta }u(\beta \log t ,x ) >1 \Big) \\&\leq e^{n(1-\delta)} \left[ \exp\left( -\frac{e^{n(\delta-\epsilon)}}{2}e^{-ct^{4+\epsilon'}}   \right)+  C^n e^{n(\delta -\epsilon)-n^3 t +n(\log 2+\gamma t) }\right].
\end{align*} Note that the summation of the right hand side of the above inequality over $n$ is finite. By the Borel-Cantelli lemma, the following holds 
 \begin{equation}\label{eq:ASevent}
  	 \mathfrak{S}^{(t)}_{\beta, \gamma} \cap  \I_{j,n}(\delta) \neq\varnothing \text{ for all } j \in [0, e^{n(1-\delta)}) \cap \Z, 
  \end{equation} with probability $1$ for all sufficiently large $n\geq 2$. Like as in the proof of Theorem~\ref{thm:Main}(a), \eqref{eq:ASevent} 
shows that $\mu_{n,1-\delta}(\mathcal{S}_{\beta}(\mathcal{V}(\gamma)))$ eventually exceeds $e^{n(1-\delta)}$ with probability $1$. From Proposition~\ref{prop:DensityTheorem}, it now follows that $\nu_{n,1-\gamma}(\mathcal{S}_{\beta}(\mathcal{V}(\gamma))) \geq C$ with probability $1$ implying that    
  	\[\Dim(\mathfrak{S}^{(t)}_{\beta,\gamma}) \geq 1-\delta, \quad \text{a.s.}
  \] Letting $\delta$ to $0$ in the above display completes the proof.

\subsubsection{Proof of Proposition~\ref{prop:inf}} \label{subsec:proof oof prop:inf}

We need to bound the probability of $u(t,x)$ getting smaller than $e^{-\gamma t}$ as $(t,x)$ varies in the set $(a,a+l_1]\times (b,b+l_2]$. The first step to get such an upper bound is to localize $u(t,x)$ in smaller boxes of bounded size. In each of the local box, one can control the probability of $u(t,x)$ taking jumps using the available tail probability of the one point distribution and the local modulus of continuity of the parabolic Anderson equation. All of these upper bounds on the probabilities of the events restricted to the smaller box will be finally combined to complete the proof. Although our proof is motivated from \cite[Proposition~3.14]{KKX18}, there is striking difference between our techniques for the valleys and those used in \cite{KKX18} for the peaks. In the latter case, there was apriori knowledge on the tail bounds of the supremum of $u(t,x)$ in an interval along spatial direction (see \cite[Theorem~5.1]{Chen16}, \cite[Proposition~3.14]{KKX18}). However, deriving such tail bound for the infimum of $u(t,x)$ is far more challenging. We will bypass this difficulty by suitable use of one point lower tail of $u(t,x)$ with its spatio-temporal modulus of continuity.


Fix $\delta>0$ and $\epsilon\in (\delta,1)$. We fix some large $a$ and $b$ which will be specified later. Let $m_1:=m_1(a,\epsilon)$ and $m_2:=m_2(a,\epsilon)$ be the least integers satisfying $m_1 e^{-a^{2-\epsilon}} > l_1$ and $m_2 e^{-a^{2-\epsilon}} > l_2$. Define   
\begin{align*}
  a_i : = a+ie^{-a^{2-\epsilon}}, \quad\text{and} \quad b_j:= b+ je^{-a^{2-\epsilon}}, \qquad I^a_i : = (a_i, a_{i+1}], \quad \text{and} \quad I^b_{j} := (b_j, b_{j+1}]
\end{align*} for any integers $i \in [0, m_1-1]$, $j\in [0,m_2-1]$ where $a_{m_1}:=a+1$ and $b_{m_2}:=b+1.$ Note that $m_1+1\leq 2l_1e^{a^{2-\epsilon}}$ and $m_2+1\leq 2l_2e^{a^{2-\epsilon}}$. 
  By the union bound,  
    \begin{align}
          \P &\Big(\text{For some }(t,x) \in (a,a+1] \times (b,b+1], u(t,x) <e^{-\gamma t} \Big)\leq (\mathbf{I})+ (\mathbf{II}), \label{eq:conti}\\ (\mathbf{I}) &:= \P \Big(\bigcup_{(i,j) \in \Z^2 \cap ([0,m_1] \times [0,m_2])}\big\{ u(a_i,b_j) <2e^{-\gamma a_i}\big\} \Big)  \nonumber\\          
       (\mathbf{II}) &:= \P \Big(\bigcup_{(i,j) \in \Z^2 \cap ([0,m_1] \times [0,m_2])}\big\{ \sup_{t,s\in I^a_i} \sup_{x,y \in I^b_j} |u(t,x)-u(s,y) | >e^{-\gamma a}\big\}\Big). \nonumber
    \end{align} 
    
    The proof will be completed by showing that $(\mathbf{I})$ and $(\mathbf{II})$ are both bounded above by $Cl_1l_2 e^{-Ca^{2-\epsilon}}$ for some constant $C>0$. The first term $(\mathbf{I})$ will be bounded above using Proposition~\ref{prop:lowertail} and the second term $(\mathbf{II})$ will be bounded using the tail bound on the spatio-temporal modulus of continuity of $u(t,x)$. We provide the details below. Applying the union bound for $\tilde{\gamma} \in (1/24,\gamma)$, we have 
    \begin{align*}
        (\mathbf{I})&\leq (m_1+1)(m_2+1) \cdot \sup_{t \in (a, a+1]}\sup_{x \in (b,b+1]}\P\big(u(t,x) < e^{-\tilde{\gamma} t} \big)\\
        &\leq 4 l_1 l_2 
        \exp \Big( 2a^{2-\epsilon} -c \left(\tilde{\gamma} -\frac{1}{24} \right)^{3-\frac{3\delta}{2}}a^{2-\delta} \Big)\leq Cl_1l_2\exp({-Ca^{2-\delta}}),
    \end{align*} where we have used $(m_1+1)(m_2+1)\leq 4l_1l_2\exp(2a^{2-\epsilon})$ and Proposition~\ref{prop:lowertail} to obtain the second inequality. The last inequality follows for all sufficiently large $a$ since $2-\delta > 2-\epsilon$. This shows the upper bound on $(\mathbf{I})$. 
    
    We now proceed to bound $(\mathbf{II})$. By \cite[Lemma~3.13]{KKX18}, there exists $N>0$ such that for all $k\in \R_{\geq 2}$, $q\in(0,1-(6/k)),$ and $\eta\in(q, 1-(6/k))$,
    \begin{equation}
        \E \Big[\sup_{t\neq s\in I^a_i} \sup_{x\neq y \in I^b_j} \Big|\frac{ u(t,x)-u(s,y) }{ ||t-s|^{1/4} + |x-y|^{1/2}|^q}\Big|^k\Big] \leq Ce^{Nk^3(a+1)}, 
    \end{equation} where $C$ depends only on $(k,\eta,q)$. Since the lengths of $|I_i^a|$ and $ |I_j^b|$ are less than  $\exp(-a^{2-\epsilon})$, we have   
\begin{equation}
    \E \Big[ \sup_{t,s\in I^a_i} \sup_{x,y \in I^b_j} |u(t,x)-u(s,y) | \Big] \leq C e^{Nk^3(a+1)} 2^{kq-1} \exp\Big(-\frac{ kq a^{2-\epsilon}}{4}\Big).
\end{equation} Applying the union bound and the above inequality yields  
\begin{align*}
    (\mathbf{II}) & \leq  (m_1+1)(m_2+1) \cdot C e^{Nk^3(a+1) + k\gamma a  }  \cdot 2^{kq} \exp\Big(-\frac{ kq a^{2-\epsilon}}{4}\Big)\\
    &\leq Cl_1l_2 \cdot 2^{kq}\exp \Big( Nk^3(a+1)+k\gamma a +a^{2-\epsilon}\big( 2-\frac{kq}{4}\big) \Big),
\end{align*} where the last inequality follows since $(m_1+1)(m_2+1)\leq 4l_1l_2\exp(2a^{2-\epsilon}).$ By choosing $k$ and $q$ such that $2- (kq/4)<0 $, the last line of the above display can be bounded above by $Cl_1l_2\exp(-Ca^{2-\epsilon})$ for all large $a>1$ with some constant $C>0$. This provides the upper bound to $(\mathbf{II})$. Substituting the upper bound on $(\mathbf{I})$ and $(\mathbf{II})$ into the right hand side of \eqref{eq:conti} completes the proof.

%
%
%
%

\section{Tail Estimates}\label{sec:Tail}
The main goal of this section is to prove Theorem~\ref{thm:lowertail}. Apart from this, we also provide the proofs of Proposition~\ref{prop:UpperTail} and~\ref{prop:lowertail} towards the end of this section. Before proceeding to main technical body of this section, we will recall some known facts and introduce few notations. Recall that the Cole-Hopf solution of the KPZ equation is none other than the logarithm of the solution of \eqref{eq:SHE}. When started from the delta initial measure, the logarithm of the solution of \eqref{eq:SHE} corresponds to the solution of the KPZ equation started from the narrow wedge initial data. We will denote the latter by $\mathcal{H}^{\mathbf{nw}}$ and define 
\begin{align}
\Upsilon_t(x) := \frac{\mathcal{H}^{\mathbf{nw}}(t,t^{2/3}x) + \frac{t}{24}}{t^{1/3}}. 
\end{align}

Below we recall few results on the tail probabilities of $\Upsilon_t$ from \cite{CG20b,CG20a,CGH21}.
 
\begin{proposition}[Theorem 1.1 of \cite{CG20a}]\label{prop:LowerTail} Fix $\delta\in(0,1/3)$ and $t_0>0$. Then, there exists $s_0=s_0 (t_0) >0 $ and $c_1=c_1(t_0), c_2=c_2(t_0)>0$ such that for all $s\geq s_0$ and $t\geq t_0$,

    \begin{equation}
       e^{-c_1s^3}+ e^{-c_1t^{1/3}s^{5/2}} \leq \P \big( \Upsilon_t(x)+\frac{x^2}{2} \leq -s  \big) \leq e^{-c_2 t^{1/3}s^{5/2}} + e^{-c_2 s^{3-\delta} } + e^{-c_2 s^3}.
    \end{equation}
\end{proposition}

%

\begin{proposition}[Proposition 1.10 of \cite{CG20b}]\label{prop:uppertail} For any $t_0>0$, there exist $s_0=s_0(t_0)>0$ and $c_1(t_0)>c_2(t_0)>0$ such that,for $t>t_0$, $s>s_0$ and $x\in \R$,
\begin{equation}
       e^{-c_1 s^{3/2}} \leq \P \big( \Upsilon_t(x)+\frac{x^2}{2} \geq s  \big) \leq e^{-c_2 s^{3/2}}. 
   \end{equation}   
\end{proposition}

\begin{proposition}[Proposition 4.2 of \cite{CGH21}]\label{prop:supupper} For any $t_0>0$ and $\nu\in[0,1)$, there exist $s_0=s_0(t_0, \nu) $ and $c_1=c_1(t_0,\nu)>c_2=c_2(t_0,\nu)>0$ such that, for $t\geq t_0 $ and $s>s_0$, 
\[
    e^{-c_1s^{3/2}} \leq \P \Big(\sup_{x\in \R} \big\{ \Upsilon_t(x) +\frac{\nu x^2}{2}\big\} \geq s \Big) \leq   e^{-c_2s^{3/2}}.
\]    
\end{proposition}

Now we proceed to prove Theorem~\ref{thm:lowertail}.

\subsection{Proof of Theorem~\ref{thm:lowertail}}

	We use Proposition~\ref{prop:Convolution} to prove this theorem. By the convolution formula (see \eqref{eq:Convolution}), it suffices to show 
	 \begin{align}\label{eq:ConvCons}
	 \P\Big(\int^{\infty}_{-\infty} u_0(x-y)e^{t^{1/3}\Upsilon_t(t^{-2/3}y)}dy\leq e^{-(\gamma - \frac{1}{24})t}\Big)\geq e^{-ct^{4+\epsilon}}
	 \end{align}
	 for all large $t$. Recall that there exists $K>0$ such that $u_0(z)\leq K$ for all $z\in \R$. Fix $\nu \in (0,1)$. It is straightforward to check that there exists $C=C(K,\gamma,t_0)>0$ such that 

\begin{equation}\label{eq:Reduc}	
\begin{aligned}
    \Big\{\sup_{x\in \R} \big\{  \Upsilon_t(x) + \frac{\nu x^2}{2} \big\}\leq - Ct^{2/3}  \Big\} &\subseteq \Big\{t^{2/3}\int^{\infty}_{-\infty} u_0(x-t^{2/3}y)e^{t^{1/3}\Upsilon_t(y)}dy\leq e^{-(\gamma - \frac{1}{24})t}\Big\} \\
    &= \Big\{\int^{\infty}_{-\infty} u_0(x-y)e^{t^{1/3}\Upsilon_t(t^{-2/3}y)}dy\leq e^{-(\gamma - \frac{1}{24})t}\Big\},
\end{aligned}
\end{equation}
	for all $t\geq t_0$. In the last equality, we used the change of variable $y\mapsto t^{-2/3}y$. Thanks to the above display, proving \eqref{eq:ConvCons} boils down to showing the following result which is certainly the main technical contribution of this section.

\begin{theorem}\label{prop:suplower} For $C>0$, $\epsilon>0$ and $\nu \in (0,1)$, there exists $c=c(t_0,C,\nu,\epsilon)>0$ and $t_0=t_0(C,\nu,\epsilon)>0$ such that for $t\geq t_0$,
\begin{equation}\label{eq:suplower}
    \P \Big( \sup_{x\in \R} \left( \Upsilon_t(x) + \frac{\nu x^2}{2}\right) \leq -Ct^{2/3} \Big) \geq \exp \big( -ct^{4+\epsilon}\big)
    \end{equation}
\end{theorem}

\begin{remark}\label{rem:InitExt}
In \eqref{eq:Reduc}, Theorem~\ref{thm:lowertail} is reduced to solving relatively simpler problem of finding lower bound to the lower tail probability of $\sup_{x\in \R}\Upsilon_t(x)$. Similar reduction is possible when the function $f:\R\to \R$ defined by $f(x):=t^{-1/3}\log u_0(t^{2/3}x)$ satisfies similar conditions as in Definition~1.1 of \cite{CG20b}. This implies Theorem~\ref{thm:lowertail} can be extended to a large class of initial data since the rest of steps in our proof does not need any specification of $u_0$.    
\end{remark}

For proving Theorem~\ref{prop:suplower}, we need the following three propositions. The first proposition will show that the supremum of $\Upsilon_t(x)+ \frac{\nu x^2}{2}$ in $x$-variable is attained with high probability in a compact interval. The second proposition will find a tail bound on the fluctuations $\Upsilon_t(x)+ \frac{\nu x^2}{2}$ on a compact interval and the last proposition is a generalization of the FKG type inequality for the KPZ equation found in \cite{CQ13}. We first state those propositions. This will be followed by the proof of Theorem~\ref{prop:suplower}. Thereafter, we will complete the proof of those propositions in three ensuing subsections.

\begin{proposition} \label{lem:supupper} For any $t_0 >0$ and $\nu \in (0,1)$, there exists $M_0 = M_0(t_0,\nu)>0$, $c=c(t_0,\nu)>0$ such that for $t\geq t_0$ and $M>M_0$,
    \begin{equation}\label{eq:Localize}
\P\Big(\mathrm{argsup}_{x\in \R} \big\{ \Upsilon_t(x) + \frac{\nu x^2}{2}\big\}\in [-M,M]\Big) \geq 1- e^{-cM^3}. 
    \end{equation}
   \end{proposition} 
\begin{proposition}\label{lem:supconti} Fix $a\in \R$ and $\epsilon\in(0,1)$. For any $t_0>0$, there exists $s_0=s_0(t_0)>0$, $c=c(t_0)>0$ such that for $t>t_0$ and $s>s_0$  
\begin{align}\label{eq:extendedModContinuity} 
     \P\Big(\sup_{x\in [a,a+\epsilon \sqrt{s}/16]} \left|\Upsilon_t(x) +\frac{x^2}{2} - \Upsilon_t(a) - \frac{a^2}{2} \right| \geq \sqrt{\epsilon} s \Big) \leq e^{-cs^{3/2}}.
\end{align}

\end{proposition}

\begin{proposition}\label{lem:fkg} Suppose $[a,b]$ and $[c,d]$ are disjoint intervals in $\R$. Then for all $s\in \R$ and $\nu>0$, we have 
       \begin{align}\label{eq:fkg}
\P &\Big(\sup_{x\in [a,b] \cup [c,d]} \left( \Upsilon_t(x) + \frac{ \nu x^2}{2}\right) \leq s  \Big)\nonumber \\ &\geq \P \Big( \sup_{x\in [a,b] } \left( \Upsilon_t(x) + \frac{\nu x^2}{2}\right) \leq s  \Big)
\cdot \P \Big( \sup_{x\in[c,d]} \left( \Upsilon_t(x) + \frac{\nu x^2}{2}\right) \leq s  \Big).
       \end{align}
\end{proposition}

\subsubsection{Proof of Theorem~\ref{prop:suplower}} 
    Choose and fix an arbitrary $\epsilon >0$. Let us define $s:= C t^{2/3}$ where $C$ is same as in Theorem~\ref{prop:suplower}. Note that 
    \begin{align}
    	\P\Big(  \sup_{x\in \R} \big\{ \Upsilon_t(x) + \frac{\nu x^2}{2} \big\}\leq -s\Big)\geq &(\mathbf{A}) - (\mathbf{B})\label{eq:cptsup},
 \end{align}
 where
 \begin{align*}    	
    	 (\mathbf{A}) &:= \P \Big( \sup_{x\in [-s^{2+\epsilon}, s^{2+\epsilon}]} \big\{ \Upsilon_t(x) + \frac{\nu x^2}{2} \big\}\leq -s\Big),\\ (\mathbf{B})&:=\P\Big(  \mathrm{argsup}_{x\in \R} \big\{\Upsilon_t(x) + \frac{\nu x^2}{2} \big\} \notin [-s^{2+\epsilon}, s^{2+\epsilon}]\Big).
    \end{align*} 
    By Proposition~\ref{lem:supupper}, we may bound $(\mathbf{B})$ from above by $\exp(-cs^{6+3\epsilon})$ whenever $s\geq s_0$ and $t\geq t_0$ where $s_0>$ and $c>0$ will only depend on some fixed $t_0>0$. When $t\geq t_0$ is sufficiently large such that $s(=Ct^{2/3})$ exceeds $s_0$, then, we can certainly write 
    \begin{align*}
    (\mathbf{B}) \leq \exp(-cs^{6+3\epsilon}). 
\end{align*}     
    Now we proceed to find lower bound to $(\mathbf{A})$. To this end, we decompose interval $[-s^{2+\epsilon},s^{2+\epsilon}]$ into smaller intervals $I_i:= [a_i,a_{i+1}]$ where 
    \[ 
    	a_i : = -s^{2+\epsilon} +(i-1)s^{-1-\epsilon} \quad \text{for } i=1,...,k-1 \quad \text{and} \quad a_k:= s^{2+\epsilon},
    \] with $k :=\lceil 32s^{3 +2\epsilon}\rceil$. Applying the FKG inequality of Proposition~\ref{lem:fkg} shows that 
    \begin{align}\label{eq:fkgI_i}
    	\P \Big( \sup_{x\in [-s^{2+\epsilon}, s^{2+\epsilon}]} \big\{ \Upsilon_t(x) + \frac{\nu x^2}{2} \big\}\leq -s\Big)\geq \prod_{i=1}^{k} \P \Big( \sup_{x\in I_i } \big\{ \Upsilon_t(x) + \frac{\nu x^2}{2} \big\}\leq -s\Big),
    \end{align} 
    Note that $|I_i| \leq s^{-1-\epsilon}/16$ for each $i$. On each interval $I_i$, by the union bound, we obtain  
    \begin{align}
     	\P \Big( \sup_{x\in I_i } \big\{ \Upsilon_t(x) + \frac{\nu x^2}{2} \big\}\leq -s\Big) \geq &\P\Big( \Upsilon_t(a_i) + \frac{a_i^2}{2} \leq -2s\Big) \nonumber\\&- \P \Big( \sup_{x\in I_i } \big|\Upsilon_t(x) +\frac{x^2}{2} - \Upsilon_t(a_i) -\frac{a_i^2}{2} \big|\geq s \Big)\label{eq:Divide}
     	\end{align} 
     	We may use Proposition~\ref{prop:LowerTail} to lower bound the first term in the right hand side of the above inequality. To this end, we get for large $s\geq s_0$ and $t\geq t_0$,  
     	\begin{align*}
     	\P\Big( \Upsilon_t(a_i) + \frac{a_i^2}{2} \leq -2s\Big)\geq \exp(- c_1t^{\frac{1}{3}}s^{\frac{5}{2}}) + \exp(-c_2 s^3). 
\end{align*}  
Since $s=Ct^{2/3}$, the right hand side is bounded below by $\exp(-cs^{3}) $ for all large $t>0$.   
To bound the second term, we use Proposition~\ref{lem:supconti}. Letting $\varepsilon:= s^{-2-4\epsilon/3}$, $\tilde{s}:=s^{2+2\epsilon/3}$ and applying stationarity of the spatial process $\Upsilon_t(x)+x^2/2$ shows  
\begin{align*}
 \P \Big( \sup_{x\in I_i } \left|\Upsilon_t(x) +\frac{x^2}{2} - \Upsilon_t(a_i) -\frac{a_i^2}{2} \right|\geq s \Big)=  \P \Big( \sup_{|x|\leq \varepsilon \sqrt{\tilde{s}}/16} \left|\Upsilon_t(x) +\frac{x^2}{2} - \Upsilon_t(0)\right|\geq \sqrt{\varepsilon}\tilde{s} \Big)\leq e^{-cs^{3+\epsilon}}
\end{align*}    	
 where the last inequality follows by applying Proposition~\ref{lem:supconti}.  Combining the bounds on both terms on the right side of \eqref{eq:Divide}, we may write 
 \begin{align}
 \text{r.h.s. of \eqref{eq:Divide}}\geq e^{-cs^{3}}
\end{align}       	
for all large $t>0$. This provides a lower bound to each of the terms in the product of \eqref{eq:fkgI_i}. Substituting those lower bounds into \eqref{eq:fkgI_i} and recalling that $k=\lceil 32s^{3+2\epsilon}\rceil$ yields 
\begin{align*}
	(\mathbf{A})=\P \Big( \sup_{x\in [-s^{2+\epsilon}, s^{2+\epsilon}]} \big\{ \Upsilon_t(x) + \frac{\nu x^2}{2} \big\}\leq -s\Big) \geq e^{-\mathfrak{c}s^{6+2\epsilon}} 
\end{align*} 
for all large $t>0$ where $\mathfrak{c}$ is a positive constant which does not depend on $t$ or $\epsilon$.
Putting the lower bound on $(\mathbf{A})$ and the upper bound on $(\mathbf{B})$ together into \eqref{eq:cptsup} shows 
\begin{equation}
	\P\Big(  \sup_{x\in \R} \big\{ \Upsilon_t(x) + \frac{\nu x^2}{2} \big\}\leq -s\Big) \geq e^{-\mathfrak{c}s^{6+2\epsilon}} - e^{-cs^{6+3\epsilon}}\geq 2^{-1}e^{-\mathfrak{c}s^{6+2\epsilon}}
\end{equation} 
for all large $t>0$. Notice that $s^{6+2\epsilon} = C^{6+2\epsilon}t^{4+4\epsilon/3}$. This completes the proof since $\epsilon$ is an arbitrary constant.


\subsubsection{Proof of Proposition~\ref{lem:supupper}}
     Note that \eqref{eq:Localize} will be proved once we show the following bound
\[
    \P(E)\leq e^{-cM^3}, \quad\text{where } E: = \Big\{\sup_{x\in \R\setminus [-M,M]} \big\{ \Upsilon_t(x) + \frac{\nu x^2}{2}\big\}  > \Upsilon_t(0)\Big\}
\] 
Fix $\nu^{\prime}\in (\nu,1)$ and choose a constant $\kappa>0$ such that $(\nu'-\nu)>2\kappa$. By the union bound, we may write 
\begin{equation}\label{eq:PE}
  \P(E) \leq \P(\sup_{x\in \R\setminus [-M,M]} \{ \Upsilon_t(x) + \frac{\nu x^2}{2}\}\geq -\kappa M^2) + \P(\Upsilon_t(0)\leq -\kappa M^2).
\end{equation} 
By Proposition~\ref{prop:LowerTail}, we may bound $\P(\Upsilon_t(0)\leq - \kappa M^2) $ by $\exp(-cM^{5})$ for all $t\geq t_0$ where $c$ is an absolute constant which only depends $\nu, \nu', t_0$. It remains to bound the first term in the right hand side of the above display. To this end, denoting $(\nu'-\nu)/2 - \kappa =: \theta$, we write  
\[\P\Big(\sup_{x\in \R\setminus [-M,M]} \{ \Upsilon_t(x) + \frac{\nu x^2}{2}\}\geq -\kappa M^2\Big) \leq  \P\Big(\sup_{x\in \R\setminus [-M,M]} \{ \Upsilon_t(x) + \frac{\nu' x^2}{2}\}\geq \theta M^2\Big)\leq \exp(-cM^3)\]
where the last inequality follows from Proposition~\ref{prop:supupper} for all $t\geq t_0$ and $M\geq M_0$. Combining the above bounds and substituting those into \eqref{eq:PE} completes the proof. 

%

\subsubsection{Proof of Proposition~\ref{lem:supconti}}
We will prove this result using Proposition~4.1 and~4.2 from \cite{DG21}. By proper identification of the notation used in this section with those used in \cite[Proposition~4.1,4.2]{DG21}, we see that there exist  $c=c(t_0)>0$ and $s_0=s_0(t_0)$ such that for any $\epsilon\in (0,1)$,
\begin{align*}
\P\Big(\inf_{|x|\leq \epsilon \sqrt{s}}\{\Upsilon_t(x)- \Upsilon_t(0)\} \leq -\sqrt{\epsilon} s\Big)&\leq e^{-cs^{3/2}}, \\
\P\Big(\sup_{|x|\leq \epsilon \sqrt{s}/16}\{\Upsilon_t(x)- \Upsilon_t(0)\} \geq \sqrt{\epsilon} s\Big)&\leq e^{-cs^{3/2}}. 
\end{align*}   
Notice that the above inequalities holds (with possibly a different constant $c$) if we replace $\Upsilon_t(x)+\frac{x^2}{2}$. The aftermath of these substitutions can be summarized in the following   
inequality: there exists $c,s_0$ depending on $t_0$ such that for all $s\geq s_0$ and $t\geq t_0$,
\begin{align*}
\P\Big(\sup_{|x|\leq \epsilon \sqrt{s}/16}\big|\Upsilon_t(x)+\frac{x^2}{2}- \Upsilon_t(0)\big|\geq \sqrt{\epsilon} s\Big)\leq e^{-cs^{3/2}}.
\end{align*}
Now \eqref{eq:extendedModContinuity} follows from the above display by the stationarity of the spatial process of $\Upsilon_t(x)+ \frac{x^2}{2}$.

\subsubsection{Proof of Proposition~\ref{lem:fkg}}
 By using the FKG inequality for the KPZ equation (see Proposition 1 in \cite{CQ13} or, Proposition 2.7 in \cite{CGH21}, we may write that any $k,k'\in \mathds{N}$, $k<k'$, $t >0 $, $x_1,...,x_k \in [a,b]$ and $x_{k+1}, \ldots ,x_{k'}\in [c,d]$, 
\[
    \P \Big( \bigcap_{l=1}^{k'} \Big\{ \Upsilon_t(x_l) + \frac{\nu x_l^2}{2} \leq s \Big\}  \Big) \geq  \P \Big( \bigcap_{l=1}^{k'} \Big\{ \Upsilon_t(x_l) + \frac{\nu x_l^2}{2} \leq s \Big\}  \Big)\cdot  \P \Big( \bigcap_{l=k'+1}^{k} \Big\{ \Upsilon_t(x_l) + \frac{\nu x_l^2}{2} \leq s \Big\}  \Big).
\] 
Suppose $x_1,\ldots ,x_k$ be the first $k$ terms of an enumerations of the rational numbers from $[a,b]$ and similarly, $x_{k+1}, \ldots , x_{k'}$ be the first $k'-k$ terms in the rational enumeration from $[c,d]$. Letting $k\to \infty$ and $k'-k\to \infty$ shows 
\begin{align*}
    \P &\Big( \bigcap_{\substack{x\in [a,b] \cup [c,d]\\x\in \Q }} \Big\{ \Upsilon_t(x) +\frac{\nu x^2}{2}\leq s \Big\} \Big)\geq \P \Big( \bigcap_{\substack{x\in [a,b] \\x\in \Q }} \Big\{ \Upsilon_t(x) +\frac{\nu x^2}{2} \leq s \Big\}\Big\} \cdot \P \Big\{ \bigcap_{\substack{x\in [c,d] \\x\in \Q }} \Big\{ \Upsilon_t(x) +\frac{\nu x^2}{2} \leq s \Big\}\Big).
\end{align*} Since $\Upsilon_t(x)$ is almost surely continuous in $x$ for any $t$, \eqref{eq:fkg} follows immediately from the above display.

\subsection{Proof of Proposition~\ref{prop:UpperTail}}\label{sec:UpperTail}

 To prove this result, we use Theorem 1.4 of \cite{CG20a} which provides upper and lower bound to the upper tail probabilities of the Cole-Hopf solution of the KPZ equation started from a large class of initial data including any bounded initial data. Since the Cole-Hopf solution is the logarithm of the solution of \eqref{eq:SHE}, those tail probability bounds are also applicable for the parabolic Anderson equation. To apply \cite[Theorem~1.4]{CG20a}, one needs $\tilde{u}_{0,t}(z):=t^{-1/3}\log u_0(t^{2/3}z)$ to satisfy two conditions of Definition~1.1 of \cite{CG20a}. The first condition requires $\tilde{u}_{0,t}(z)$ to be upper bounded by a parabola in $z$-variable for all $t\geq t_0$ and the second condition requires $\tilde{u}_{0,t}(z)$ to be lower bounded by a finite constant in an interval around $0$. If $u_0$ is a bounded positive initial data, then, the above two conditions will be satisfied for $\tilde{u}_{0,t}(\cdot)$. Hence,  by Theorem~1.4 of \cite{CG20a}, for any $t_0>0$, there exists $c_1>c_2>0$ and $s_0$ such that for all $t>t_0$ and $s>s_0$, 
 \begin{align}\label{eq:UpLowBd}
 e^{ -c_1s^{3/2}}\leq  \P\left(u(t,0)>e^{-\frac{t}{24}+s\left(\frac{t}{2}\right)^{1/3}+\frac{2}{3}\log t  } \right)\leq e^{ -c_2s^{3/2}}.
\end{align}     
From the above inequalities, \eqref{eq:lya} for $x=0$ follows by substituting $s=2^{1/3} \left((\frac{1}{24}-\gamma)t^{2/3} - \frac{2\log t}{3t^{1/3}} \right)$. Since $u_0$ is a bounded positive initial data, the inequalities in \eqref{eq:UpLowBd} also hold for $x\neq 0$ where the corresponding constants $c_1,c_2$ would not depend on $x$. The reason behind this fact lies in the convolution formula \eqref{eq:Convolution} for the one point distribution of $u(t,x)$ and the proof can be completed using similar argument as in \cite[Theorem~1.4]{CG20a}. Once \eqref{eq:UpLowBd} holds for $x\neq 0$, we obtain \eqref{eq:lya} for $x\neq 0$ by making the appropriate substitution of $s$ by $\gamma$ as prescribed above. This completes the proof. 

\subsection{Proof of Proposition~\ref{prop:lowertail}}\label{sec:LowerTail} 
 The proof of this proposition is analogous to how we have proved Theorem~\ref{prop:UpperTail} using Theorem~1.4 of \cite{CG20a}. In this case, we need to use Theorem~1.2 of \cite{CG20a}. The arguments are similar as above. We skip the details for brevity.

\section*{Acknowledgments} The first author thanks Sayan Das for many useful conversations and numerous suggestions on the first draft of this paper. 
The second author thanks Professor Kunwoo Kim for valuable discussions with his continuous support and encouragement.
The second author was supported by the NRF (National Research Foundation of Korea) grants 2019R1A5A1028324 and 2020R1A2C4002077.

\bibliographystyle{alpha}		
\bibliography{refs_valley}

%
%
%
%
%
%

\end{document}